\documentclass{article}
\usepackage{amsfonts}
\usepackage{amssymb}
\usepackage{amsmath}

\input{tcilatex}

\begin{document}

\begin{center}
\textbf{THE INDEPENDENCE OF }$p$

\textbf{OF THE LIPSCOMB'S }$L(A)$\textbf{\ SPACE FRACTALIZED IN }$l^{p}(A)$

\bigskip

Radu MICULESCU and Alexandru MIHAIL
\end{center}

\bigskip

\textbf{ABSTRACT.} In one of our previous papers we proved that, for an
infinite set $A$ and $p\in \lbrack 1,\infty )$, the embedded version of the
Lipscomb's space $L(A)$ in $l^{p}(A)$, $p\in \lbrack 1,\infty )$, with the
metric induced from $l^{p}(A)$, denoted by $\omega _{p}^{A}$, is the
attractor of an infinite iterated function system comprising affine
transformations of $l^{p}(A)$. In the present paper we point out that $%
\omega _{p}^{A}=\omega _{q}^{A}$, for all $p,q\in \lbrack 1,\infty )$ and,
by providing a complete description of the convergent sequences from $\omega
_{p}^{A}$, we prove that the topological structure of $\omega _{p}^{A}$\ is
independent of $p$.

2000\ Mathematics Subject Classification: 37C70, 28A80, 54A20, 54B15

Key words and phrases:\textit{\ }{\small Lipscomb's space, convergent
sequences,}\textit{\ }{\small infinite iterated function system}

\bigskip

\begin{center}
1. INTRODUCTION

\bigskip
\end{center}

\bigskip

N\"{o}beling's classical imbedding theorem (see [3] and [7]) states that a $%
n $-dimensional separable metric space can be topologically imbedded in the
product $I^{2n+1}$ of $2n+1$ copies of the one-dimensional unit interval $I$%
. Lipscomb's space $L(A)$ was introduced in order to solve the long standing
problem of finding an analogue to N\"{o}beling's theorem. More precisely,
all finite-dimensional metric spaces of weight $\left\vert A\right\vert \geq
\aleph _{0}$ are modeled as subspaces of finite products of $L(A)$ (see [2]
and [3]).

On one hand, J. C. Perry (see [8]) embedded $L(A)$ in the Tikhonov cube $%
I^{A}$ and showed that this set with the topology induced by $I^{A}$ is the
attractor of an iterated function system containing an infinite number of
affine transformations of $I^{A}$. In this way, $L(A)$ arguably provided the
first notable example where an infinite IFS plays a key role.

On the other hand, S. L. Lipscomb and J. C. Perry\ (see [1]) showed that $%
L(A)$ can be embedded in $l^{2}(A)$. R. Miculescu and A. Mihail (see [4])
showed that the imbedded version of $L(A)$ endowed with the $l^{2}(A)$%
-induced topology is the attractor of an infinite iterated function system
comprising affine transformations of $l^{2}(A)$.

In [6], by using some results concerning the shift space for an infinite IFS
(see [5]), we proved that the embedded version of $L(A)$ in $l^{p}(A)$, $%
p\in \lbrack 1,\infty )$, with the metric induced from $l^{p}(A)$, denoted
by $\omega _{p}^{A}$, is the attractor of an infinite iterated function
system comprising affine transformations of $l^{p}(A)$.

In the present paper we give a complete description of the convergent
sequences from $\omega _{p}^{A}$. As it turns out that $\omega
_{p}^{A}=\omega _{q}^{A}$, for all $p,q\in \lbrack 1,\infty )$, by using the
above mentioned description we infer that the topological structure of $%
\omega _{p}^{A}$\ is independent of $p$.

\bigskip

\begin{center}
2. PRELIMINARIES
\end{center}

\bigskip

\textbf{Definition 2.1}. \textit{Let }$(X,d_{X})$\textit{\ and }$(Y,d_{Y})$%
\textit{\ be two metric spaces. A family of functions }$(f_{i})_{i\in I}$%
\textit{\ from }$X$\textit{\ to }$Y$\textit{\ is called bounded if the set }$%
\bigcup\limits_{i\in I}f_{i}(A)$\textit{\ is bounded, for every bounded
subset }$A$\textit{\ of }$X$\textit{.}

\bigskip

\textbf{Theorem 2.1}. \textit{If }$(X,d)$\textit{\ is a complete metric
space, then }$(\mathcal{B}^{\ast }(X),h)$\textit{\ is a complete metric
space, where }$h$\textit{\ is the Hausdorff-Pompeiu distance and }$\mathcal{B%
}^{\ast }(X)$\textit{\ is the set of non-empty bounded closed subsets of }$X$%
\textit{.}

\bigskip

\textbf{Definition 2.2}. \textit{An infinite iterated function system (IIFS)
on }$X$\textit{\ consists of a bounded family of contractions }$%
(f_{i})_{i\in I}$\textit{\ on }$X$\textit{\ such that }$\underset{i\in I}{%
\sup }Lip(f_{i})<1$\textit{\ and it is denoted by} $\mathcal{S}%
=(X,(f_{i})_{i\in I})$\textit{.}

\bigskip

One can associate to an infinite iterated function system $\mathcal{S}%
=(X,(f_{i})_{i\in I})$ the function $F_{\mathcal{S}}:\mathcal{B}^{\ast
}(X)\rightarrow \mathcal{B}^{\ast }(X)$ given by 
\begin{equation*}
F_{\mathcal{S}}(B)=\overline{\bigcup\limits_{i\in I}f_{i}(B)}\text{,}
\end{equation*}
for all $B\in \mathcal{B}^{\ast }(X)$\textit{.}

Therefore, since $Lip\mathbf{(}F_{\mathcal{S}})\leq \underset{i\in I}{\sup }%
Lip(f_{i})$, $F_{\mathcal{S}}$ is a contraction and using the Banach's
contraction theorem, one can prove the following:

\bigskip

\textbf{Theorem 2.2.} \textit{Given a complete metric space} $(X,d)$ \textit{%
and an IIFS} $\mathcal{S}=(X,(f_{i})_{i\in I})$\textit{, there exists a
unique} $A(\mathcal{S})\in \mathcal{B}^{\ast }(X)$ \textit{such that} 
\begin{equation*}
F_{\mathcal{S}}(A(\mathcal{S}))=A(\mathcal{S})\text{.}
\end{equation*}

\textit{The set }$A(\mathcal{S})$ \textit{is called the attractor associated
to }$\mathcal{S}$.

\bigskip

In the following $\mathbb{N}$ denotes the natural numbers and $\mathbb{N}%
^{\ast }\mathbb{=N}-\{0\}$. Given two sets $A$ and $B$, by $B^{A}$ we mean
the set of functions from $A$ to $B$. By $\Lambda =\Lambda (B)$ we mean the
set $B^{\mathbb{N}^{\ast }}$ and by $\Lambda _{n}=\Lambda _{n}(B)$ we mean
the set $B^{\{1,2,...,n\}}$. The elements of $\Lambda =\Lambda (B)=B^{%
\mathbb{N}^{\ast }}$ are written as infinite words $\omega =\omega
_{1}\omega _{2}...\omega _{m}\omega _{m+1}...$ and the elements of $\Lambda
_{n}=\Lambda _{n}(B)=B^{\{1,2,...,n\}}$ are written as words $\omega =\omega
_{1}\omega _{2}...\omega _{n}$. Hence $\Lambda (B)$ is the set of infinite
words with letters from the alphabet $B$ and $\Lambda _{n}(B)\,$is the set
of words of length $n$ with letters from the alphabet $B$. By $\Lambda
^{\ast }=\Lambda ^{\ast }(B)$ we denote the set of all finite words $\Lambda
^{\ast }=\Lambda ^{\ast }(B)=\bigcup\limits_{n\in \mathbb{N}^{\ast }}\Lambda
_{n}(B)\cup \{\lambda \}$, where by $\lambda $ we mean the empty word. If $%
\omega =\omega _{1}\omega _{2}...\omega _{m}\omega _{m+1}...\in \Lambda (B)$
or if $\omega =\omega _{1}\omega _{2}...\omega _{n}\in \Lambda _{n}(B)$,
where $m,n\in \mathbb{N}^{\ast }$, $n\geq m$, then the word $\omega
_{1}\omega _{2}...\omega _{m}$ is denoted by $[\omega ]_{m}$.

For a nonvoid set $I$, on $\Lambda =\Lambda (I)=(I)^{\mathbb{N}^{\ast }}$,
we consider the metric 
\begin{equation*}
d_{\Lambda }(\alpha ,\beta )=\overset{\infty }{\underset{k=1}{\sum }}\frac{%
1-\delta _{\alpha _{k}}^{\beta _{k}}}{3^{k}}\text{,}
\end{equation*}%
where 
\begin{equation*}
\delta _{x}^{y}=\left\{ 
\begin{array}{c}
1\text{, if }x=y \\ 
0\text{, if }x\neq y%
\end{array}%
\right. \text{.}
\end{equation*}

\bigskip

Let $(X,d)$ be a metric space, $\mathcal{S}=(X,(f_{i})_{i\in I})$ be an IIFS
on $X$ and $A\overset{not}{=}A(\mathcal{S})$ its attractor. For $\omega
=\omega _{1}\omega _{2}...\omega _{m}\in \Lambda _{m}(I)$, we consider $%
f_{\omega }\overset{not}{=}f_{\omega _{1}}\circ f_{\omega _{2}}\circ \ldots
\circ f_{\omega _{m}}$ and, for a subset $H$ of $X$, $H_{\omega }\overset{not%
}{=}f_{\omega }(H)$. In particular $A_{\omega }=f_{\omega }(A)$.

We also consider $f_{\lambda }=Id$ and $A_{\lambda }=A$.

\bigskip

The following result is part of Theorem 4.1 from [5].

\bigskip

\textbf{Theorem 2.3.} \textit{Then,} \textit{for every} $a\in A$ \textit{and
every} $\omega \in \Lambda $\textit{, the set }$\bigcap\limits_{m\in \mathbb{%
N}^{\ast }}\overline{A_{[\omega ]_{m}}}$ \textit{consists on a single
element denoted by }$a_{\omega }$\textit{\ and }%
\begin{equation*}
\underset{m\rightarrow \infty }{\lim }f_{[\omega ]_{m}}(a)=a_{\omega }\text{.%
}
\end{equation*}

\textit{In the framework of the above theorem,} \textit{the function} $\pi
:\Lambda \rightarrow A$\textit{, defined by} 
\begin{equation*}
\pi (\omega )=a_{\omega }\text{,}
\end{equation*}%
\textit{for every} $\omega \in \Lambda $\textit{\ (which is called the
canonical projection from the shift space on the attractor of the IIFS) is
continuos.}

\bigskip

Let $A$ be an arbitrary infinite set and single out a point $z$ of $A$. Let
us consider the set $A^{^{\prime }}=A-\{z\}$. For a given $p\in \lbrack
1,\infty )$, the points of $l^{p}(A)$ are collections of real numbers
indexed by points of $A^{^{\prime }}$. If $E$ is the set of real numbers,
then $x\in l^{p}(A)$ means $x=\{x_{a}\}\in E^{A^{^{\prime }}}$ such that $%
x_{a}=0$ for all but countable many $a\in A^{^{\prime }}$ and $%
\sum\limits_{a}\left\vert x_{a}\right\vert ^{p}$ converges. The topology of $%
l^{p}(A)$ is induced from the metric $d_{p}(x,y)=\left(
\sum\limits_{a}\left\vert x_{a}-y_{a}\right\vert ^{p}\right) ^{\frac{1}{p}}$%
, where we think $x_{a}$ as the $a$-th coordinate of $x$. By $\left\Vert
x\right\Vert _{p}$ we mean $d_{p}(x,0)$.

Let us also consider, for the case when $A$ is an arbitrary set with the
discrete topology, the Baire space $N(A)$ which is the topological product
of countably many copies $A_{n}$ of $A$. Hence the points of $N(A)$ consist
of all sequences $v=a_{1}a_{2}...a_{n}...$, with $a_{n}\in A$. Moreover $%
N(A) $ is a metric space with the metric 
\begin{equation*}
d(v\mathbf{,}v^{^{\prime }})=\left\{ 
\begin{array}{cc}
\dfrac{1}{k}, & \text{if }k\text{ is the first index where }a_{i}\neq
a_{i}^{^{\prime }} \\ 
0, & \text{if }v\mathbf{=}v^{^{\prime }}%
\end{array}%
\right. \text{.}
\end{equation*}

Let us note that $N(A)=\Lambda (A)$ and that the metrics $d$ and $d_{\Lambda
}$ are equivalent.

Lipscomb's space $L(A)$ is a quotient space of $N(A)$ such that each
equi\-va\-lence class consists of either a single point or two points. Those
classes with two points come from identifying the point $%
a_{1}a_{2}...a_{k-2}a_{k-1}a_{k}a_{k}...$ with the point $%
a_{1}a_{2}...a_{k-2}a_{k}a_{k-1}a_{k-1}...$, where $a_{k-1}\neq a_{k}$.
Therefore Lipscomb's space $L(A)$ is obtained via a projection or
identification map $p:N(A)\rightarrow L(A)$.

For $p\in \lbrack 1,\infty )$, let us consider the function $%
p_{p}:N(A)\rightarrow l^{p}(A)$ given by $p_{p}(\alpha )=(\alpha _{b})_{b\in
A^{^{\prime }}}$, where $\alpha =a_{1}a_{2}....\in N(A)$ and 
\begin{equation*}
\alpha _{b}=\left\{ 
\begin{array}{cc}
\sum\limits_{k\text{ with }a_{k}=b}\dfrac{1}{2^{k}}, & \text{if there exists 
}k\text{ such that }a_{k}=b \\ 
0, & \text{if there exists no }k\text{ such that }a_{k}=b%
\end{array}%
\right. \text{.}
\end{equation*}

Moreover, we consider the function $s_{p}:L(A)\rightarrow l^{p}(A)$ given by 
$s_{p}(\widehat{\alpha })=p_{p}(\alpha )$, for every $\widehat{\alpha }\in
L(A)$.

Let us note that $s_{p}$ \textit{is well-defined}.

Indeed, let us consider $\alpha ,\beta \in N(A)$ which are equivalent. We
want to prove that $p_{p}(\alpha )=p_{p}(\beta )$. If $\alpha =\beta $ there
is nothing to prove. If $\alpha =a_{1}a_{2}...a_{k-2}a_{k-1}a_{k}a_{k}...$
and $\beta =a_{1}a_{2}...a_{k-2}a_{k}a_{k-1}a_{k-1}...$, where $a_{k-1}\neq
a_{k}$, then:

i) 
\begin{equation*}
\alpha _{b}=\beta _{b}=0\text{,}
\end{equation*}
if $b\notin \{a_{1},a_{2},...,a_{k-2},a_{k-1},a_{k}\}$;

ii) 
\begin{equation*}
\alpha _{a_{k-1}}=\sum\limits_{i\text{ with }a_{i}=a_{k-1}}\dfrac{1}{2^{i}}%
=\sum\limits_{i\in \{1,2,...,k-2\}\text{ with }a_{i}=a_{k-1}}\dfrac{1}{2^{i}}%
+\dfrac{1}{2^{k-1}}=
\end{equation*}%
\begin{equation*}
=\sum\limits_{i\in \{1,2,...,k-2\}\text{ with }a_{i}=a_{k-1}}\dfrac{1}{2^{i}}%
+\overset{\infty }{\underset{i=k+1}{\sum }}\dfrac{1}{2^{i}}=\beta _{a_{k-1}}%
\text{;}
\end{equation*}

iii) 
\begin{equation*}
\alpha _{a_{k}}=\sum\limits_{i\text{ with }a_{i}=a_{k}}\dfrac{1}{2^{i}}%
=\sum\limits_{i\in \{1,2,...,k-2\}\text{ with }a_{i}=a_{k}}\dfrac{1}{2^{i}}+%
\overset{\infty }{\underset{i=k}{\sum }}\dfrac{1}{2^{i}}=
\end{equation*}%
\begin{equation*}
=\sum\limits_{i\in \{1,2,...,k-2\}\text{ with }a_{i}=a_{k}}\dfrac{1}{2^{i}}+%
\dfrac{1}{2^{k-1}}=\beta _{a_{k}}\text{;}
\end{equation*}

iv) 
\begin{equation*}
\alpha _{b}=\beta _{b}=\sum\limits_{i\in \{1,2,...,k-2\}\text{ with }a_{i}=b}%
\dfrac{1}{2^{i}}\text{,}
\end{equation*}%
if $b\in \{a_{1},a_{2},...,a_{k-2}\}-\{a_{k-1},a_{k}\}$.

Hence, $s_{p}$ is well-defined.

\bigskip

\textbf{Proposition 2.1}. $s_{p}$ \textit{is injective.}

\textit{Proof}. Let us consider $\alpha =a_{1}a_{2}...a_{n}...$ and $\beta
=b_{1}b_{2}...b_{n}...$ two arbitrary elements of $N(A)$ such that 
\begin{equation*}
s_{p}(\widehat{\alpha })=s_{p}(\widehat{\beta })\text{,}
\end{equation*}
i.e. 
\begin{equation*}
\alpha _{b}=\beta _{b}\text{,}
\end{equation*}
for each $b\in A^{^{\prime }}$.

We want to prove that $\alpha $ and $\beta $ are equivalent.

Let us suppose that $a_{1}\neq b_{1}$. In this case, let us note that if $%
a_{1}\neq z$, then taking into account that $\alpha _{a_{1}}=\beta _{a_{1}}$%
, we get $\frac{1}{2}\leq \sum\limits_{k\text{ with }a_{k}=a_{1}}\dfrac{1}{%
2^{k}}=\sum\limits_{k\text{ with }b_{k}=a_{1}}\dfrac{1}{2^{k}}%
=\sum\limits_{k>1\text{ with }b_{k}=a_{1}}\dfrac{1}{2^{k}}\leq
\sum\limits_{k>1}\dfrac{1}{2^{k}}=\frac{1}{2}$, which implies that 
\begin{equation}
a_{k}\neq a_{1}\text{ and }b_{k}=a_{1}\text{,}  \tag{*}
\end{equation}%
for all $k>1$. In a similar manner, we obtain that if $b_{1}\neq z$, then%
\begin{equation}
b_{k}\neq b_{1}\text{ and }a_{k}=b_{1}\text{,}  \tag{**}
\end{equation}%
for all $k>1$. If $z\notin \{a_{1},b_{1}\}$, then, according to $(\ast \ast
) $ and $(\ast )$, we get $\alpha =a_{1}b_{1}...b_{1}...$ and $\beta
=b_{1}a_{1}...a_{1}...$, which assure us that $\alpha $ and $\beta $ are
equivalent. If $z\in \{a_{1},b_{1}\}$, we can suppose, without loss of
generality, that $a_{1}=z\neq b_{1}$. Then, according to $(\ast \ast )$, we
get $\alpha =zb_{1}...b_{1}...$, and $b_{k}\neq b_{1}$, for all $k>1$. We
claim that $b_{k}=z$, for all $k>1$. Indeed, if there exists $k_{0}>1$ such
that $b_{k_{0}}\neq z$, the we get the following contradiction: $0=\alpha
_{b_{k_{0}}}=\beta _{b_{k_{0}}}\geq \frac{1}{2^{_{k_{0}}}}$. Hence $\beta
=b_{1}zz...z...$, and again we infer that $\alpha $ and $\beta $ are
equivalent.

Consequently, if $a_{1}\neq b_{1}$, then $\alpha $ and $\beta $ are
equivalent.

Therefore we can suppose that $\alpha =a_{1}a_{2}...a_{n}...$ and $\beta
=a_{1}b_{2}...b_{n}...$ .

Let us suppose that $a_{2}\neq b_{2}$. In this case, let us note that if $%
a_{2}\neq z$, then taking into account that $\alpha _{a_{2}}=\beta _{a_{2}}$%
, we get $\sum\limits_{k\text{ with }a_{k}=a_{2}}\dfrac{1}{2^{k}}%
=\sum\limits_{k\text{ with }b_{k}=a_{2}}\dfrac{1}{2^{k}}$ and therefore%
\begin{equation*}
\frac{1}{2^{2}}\leq \sum\limits_{k>1\text{ with }a_{k}=a_{2}}\dfrac{1}{2^{k}}%
=\sum\limits_{k>1\text{ with }b_{k}=a_{2}}\dfrac{1}{2^{k}}=\sum\limits_{k>2%
\text{ with }b_{k}=a_{2}}\dfrac{1}{2^{k}}\leq \sum\limits_{k>2}\dfrac{1}{%
2^{k}}=\frac{1}{2^{2}}\text{,}
\end{equation*}%
which implies that%
\begin{equation}
a_{k}\neq a_{2}\text{ and }b_{k}=a_{2}\text{,}  \tag{***}
\end{equation}%
for all $k>2$. In a similar manner, we obtain that if $b_{2}\neq z$, then%
\begin{equation}
b_{k}\neq b_{2}\text{ and }a_{k}=b_{2}\text{,}  \tag{****}
\end{equation}%
for all $k>2$. If $z\notin \{a_{2},b_{2}\}$, then, according to $(\ast \ast
\ast \ast )$ and $(\ast \ast \ast )$, we get $\alpha
=a_{1}a_{2}b_{2}b_{2}...b_{2}...$ and $\beta
=a_{1}b_{2}a_{2}a_{2}...a_{2}... $, which assure us that $\alpha $ and $%
\beta $ are equivalent. If $z\in \{a_{2},b_{2}\}$, we can suppose, without
loss of generality, that $a_{2}=z\neq b_{2}$. Then, according to $(\ast \ast
\ast \ast )$, we get $\alpha =a_{1}zb_{2}b_{2}...b_{2}...$, and $b_{k}\neq
b_{2}$, for all $k>2$. We claim that $b_{k}=z$, for all $k>2$. Indeed, if
there exists $k_{0}>2$ such that $b_{k_{0}}\neq z$, the we get the following
contradiction: $\frac{1}{2}\geq \alpha _{b_{k_{0}}}=\beta _{b_{k_{0}}}\geq 
\frac{1}{2^{_{k_{0}}}}$. Hence $\beta =a_{1}b_{2}zz...z...$, and again we
infer that $\alpha $ and $\beta $ are equivalent.

Consequently, if $a_{2}\neq b_{2}$, then $\alpha $ and $\beta $ are
equivalent.

Therefore we can suppose that $\alpha =a_{1}a_{2}a_{3}...a_{n}...$ and $%
\beta =a_{1}a_{2}b_{3}...b_{n}...$ .

Continuing this procedure, we obtain that $\alpha $ and $\beta $ are
equivalent. $\square $

\bigskip

The set $p_{p}(N(A))=s_{p}(L(A))$ is denoted by $\omega _{p}^{A}$.

For each $a\in A$, let $f_{a}:l^{p}(A)\rightarrow l^{p}(A)$ be the function
given by 
\begin{equation*}
f_{a}(x)=\frac{1}{2}(x+u_{a}),
\end{equation*}
for all $x\in l^{p}(A)$, where $u_{z}=0_{l^{p}(A)}\in \Delta _{p}^{A}$ and,
for $a\in A\setminus \{z\}$, $u_{a}=(\alpha _{j})_{j\in A^{^{\prime }}}\in
\Delta _{p}^{A}$ is described by $\alpha _{j}=0$, for $j\neq a$ and $\alpha
_{a}=1$.

Let us consider the IIFS $\mathcal{S}=(l^{p}(A),(f_{a})_{a\in A})$.
According to Theorem 2.2, there exists a bounded closed non-empty subset $M$
of $l^{p}(A)$, called the attractor of $\mathcal{S}$, such that 
\begin{equation*}
M=F_{\mathcal{S}}(M)=\overline{\bigcup\limits_{a\in A}f_{a}(M)}\text{.}
\end{equation*}

Let us note that $N(A)$ can be seen as the shift space for the IIFS $%
\mathcal{S}$.

\bigskip

The following result is Theorem 6.1 from [6].

\bigskip

\textbf{Theorem 2.4. }\textit{With the above notations }$\omega
_{p}^{A}=p_{p}(N(A))=s_{p}(L(A))$\textit{\ is the attractor of the IIFS }$%
\mathcal{S}=(l^{p}(A),(f_{a})_{a\in A})$\textit{\ and }$p_{p}:N(A)%
\rightarrow \omega _{p}^{A}$\textit{\ is the canonical projection from the
shift space on the attractor of the IIFS\ }$\mathcal{S}$\textit{.}

\bigskip

For more details concerning infinite iterated function systems (IIFSs) and
the shift space associated to an IIFS one can consult [5]. Likewise, more
details about $\omega _{p}^{A}$ as the attractor of an infinite iterated
function system comprising affine transformations of $l^{p}(A)$ can be found
in [6].

\bigskip

\begin{center}
3. THE\ MAIN\ RESULTS

\bigskip
\end{center}

\textbf{Proposition 3.1}. \textit{The equality} 
\begin{equation*}
\omega _{p}^{A}=\omega _{q}^{A}
\end{equation*}%
\textit{is valid for all }$p,q\in \lbrack 1,\infty )$.

\textit{Proof}. According to the proof of Theorem 2.4,%
\begin{equation*}
0_{l^{p}(A)}\in \omega _{p}^{A}\text{.}
\end{equation*}
\ and $p_{p}=\pi $, which implies%
\begin{equation*}
\omega _{p}^{A}=p_{p}(N(A))=\pi (N(A))\text{.}
\end{equation*}

Taking into account Theorem 2.3, we have%
\begin{equation*}
\pi (\omega )=\underset{m\longrightarrow \infty }{\lim }f_{[\omega
]_{m}}(0_{l^{p}(A)})\text{,}
\end{equation*}%
for each $\omega \in N(A)$, so%
\begin{equation*}
\omega _{p}^{A}=\{\underset{m\longrightarrow \infty }{\lim }f_{[\omega
]_{m}}(0_{l^{p}(A)})\mid \omega \in N(A)\}\text{.}
\end{equation*}

As all the functions $f_{a}$ from the definition of the IIFS whose attractor
is $\omega _{p}^{A}$ do not depend on $p$, we conclude that 
\begin{equation*}
\omega _{p}^{A}=\omega _{q}^{A}\text{,}
\end{equation*}%
for all $p,q\in \lbrack 1,\infty )$. $\square $

\bigskip

\textbf{Proposition 3.2}. \textit{The function }$p_{p}:N(A)\rightarrow
l^{p}(A)$\textit{\ is continuous.}

\textit{Proof}. It results from the fact that $p_{p}$ is $\pi $ and one can
use Theorem 2.3. $\square $

\bigskip

\textbf{Proposition 3.3}. \textit{Let }$(x_{n})_{n\in \mathbb{N}^{\ast }}$ 
\textit{be a sequence of elements from }$\omega _{p}^{A}$\textit{, where} $%
x_{n}=p_{p}(\alpha ^{n})=(\alpha _{a}^{n})_{a\in A^{\prime }}$\textit{,} 
\textit{with }$\alpha ^{n}=a_{1}^{n}a_{2}^{n}...\in N(A)$,\textit{\ and }$%
x=p_{p}(\alpha )=(\alpha _{a})_{a\in A^{\prime }}\in \omega _{p}^{A}$\textit{%
, where }$\alpha =a_{1}a_{2}...\in N(A)$\textit{.}

\textit{If for every }$q\in \mathbb{N}^{\ast }$ \textit{the sequence }$%
(a_{n})_{n\in \mathbb{N}^{\ast }}$\textit{\ is not constant after the rank }$%
q$ \textit{and}

\begin{equation*}
\underset{n\rightarrow \infty }{\lim }||x_{n}-x||_{p}=0\text{,}
\end{equation*}%
\textit{then}%
\begin{equation*}
\underset{n\rightarrow \infty }{\lim }\alpha ^{n}=\alpha \text{,}
\end{equation*}%
\textit{in }$N(A)$\textit{.}

\textit{Proof}. Let us consider, for $m\in \mathbb{N}^{\ast }$, the
following statement:

P($m$): There exists $n_{m}$ such that $a_{k}^{n}=a_{k}$ for every $n\geq
n_{m}$ and every $k\in \{1,2,...,m\}$.

We shall prove, using the method of mathematical induction, that P($m$) is
true for every $m\in \mathbb{N}^{\ast }$.

This implies that 
\begin{equation*}
\underset{n\rightarrow \infty }{\lim }\alpha ^{n}=\alpha \text{,}
\end{equation*}%
in\textit{\ }$N(A)$\textit{.}

Indeed, for every $\varepsilon >0$ let us choose $m_{\varepsilon }\in 
\mathbb{N}^{\ast }$ such that $\frac{1}{m_{\varepsilon }}<\varepsilon $.
Since $P(m_{\varepsilon })$ is true, there exists $n_{m_{\varepsilon }}%
\overset{not}{=}n_{\varepsilon }$ such that, for every $n\in \mathbb{N}%
^{\ast }$, $n\geq n_{\varepsilon }$, the first index $k$ for which $%
a_{k}^{n}\neq a_{k}$ is greater than $m_{\varepsilon }$. In other words, for
every $\varepsilon >0$ there exists $n_{\varepsilon }\in \mathbb{N}^{\ast }$
such that $d(\alpha ^{n},\alpha )\leq \frac{1}{m_{\varepsilon }}<\varepsilon 
$, for every $n\in \mathbb{N}^{\ast }$, $n\geq n_{\varepsilon }$.

Let us prove that P($1$) is true.

If this is not the case, then for every $n\in \mathbb{N}^{\ast }$ there
exists $n^{\prime }\in \mathbb{N}^{\ast }$, $n^{\prime }\geq n$, such that $%
a_{1}^{n^{\prime }}\neq a_{1}$ and therefore, by passing to a subsequence,
we can suppose that $a_{1}^{n}\neq a_{1}$, for every $n\in \mathbb{N}^{\ast
} $.

We are going to treat two cases.

In the first one, which is described by the situation that $\{n\in \mathbb{N}%
^{\ast }\mid a_{1}^{n}\neq z\}$ is infinite, by passing to a subsequence, we
can suppose that $a_{1}^{n}\neq z$, for every $n\in \mathbb{N}^{\ast }$.

Then, for all $n\in \mathbb{N}^{\ast }$, with the convention that $\underset{%
k\text{ with }a_{k}^{n}=a}{\sum }\frac{1}{2^{k}}$ (respectively $\underset{k%
\text{ with }a_{k}=a}{\sum }\frac{1}{2^{k}}$) is $0$ if there is no $k$ such
that $a_{k}^{n}=a$ (respectively $a_{k}=a$), we have

\begin{equation*}
|\alpha _{a_{1}^{n}}^{n}-\alpha _{a_{1}^{n}}|=|\sum_{k\text{ with }%
a_{k}^{n}=a_{1}^{n}}\frac{1}{2^{k}}-\sum_{k\text{ with }a_{k}=a_{1}^{n}}%
\frac{1}{2^{k}}|=
\end{equation*}

\begin{equation}
=|\sum_{k\text{ with }a_{k}^{n}=a_{1}^{n}}\frac{1}{2^{k}}-\sum_{k>1\text{
with }a_{k}=a_{1}^{n}}\frac{1}{2^{k}}|\geq \frac{1}{2}-\sum_{k>1\text{ with }%
a_{k}=a_{1}^{n}}\frac{1}{2^{k}}\text{.}  \tag{*}
\end{equation}

We claim that the sequence $(a_{n})_{n\in \mathbb{N}^{\ast }}$ is constant.
This contradicts the fact that for every $q\in \mathbb{N}^{\ast }$ this
sequence\ is not constant after the rank $q$.

The claim is true since otherwise there exist $s,t\in \mathbb{N}^{\ast }$, $%
s<t$, such that $a_{s}\neq a_{t}$. Then let us choose $l,n\in \mathbb{N}%
^{\ast }$ such that $s<t<l$ and 
\begin{equation*}
\left\Vert x_{n}-x\right\Vert _{p}<\frac{1}{2^{m}}\text{,}
\end{equation*}%
for all $m\in \{1,2,...,l\}$. If for a given $m\in \{1,2,...,l\}$ we have $%
a_{m}\neq a_{1}^{n}$, then, since%
\begin{equation*}
\sum_{k>1\text{ with }a_{k}=a_{1}^{n}}\frac{1}{2^{k}}=\sum_{k>1\text{ with }%
a_{k}=a_{1}^{n}\text{ and }k\neq m}\frac{1}{2^{k}}<1-\frac{1}{2^{m}}\text{,}
\end{equation*}%
we get, using $(\ast )$, the following contradiction%
\begin{equation*}
\frac{1}{2^{m}}<1-\sum_{k>1\text{ with }a_{k}=a_{1}^{n}}\frac{1}{2^{k}}\leq
|\alpha _{a_{1}^{n}}^{n}-\alpha _{a_{1}^{n}}|\leq \left\Vert
x_{n}-x\right\Vert _{p}<\frac{1}{2^{m}}\text{.}
\end{equation*}%
Consequently we infer that 
\begin{equation*}
a_{1}=a_{2}=...=a_{s}=...=a_{t}=...=a_{l}=a_{1}^{n}
\end{equation*}%
which contradicts the fact that 
\begin{equation*}
a_{s}\neq a_{t}\text{.}
\end{equation*}

In the second case, which is described by the situation that $\{n\in \mathbb{%
N}^{\ast }\mid a_{1}^{n}\neq z\}$ is finite, by passing to a subsequence, we
can suppose that $a_{1}^{n}=z$, for every $n\in \mathbb{N}^{\ast }$.

Then, for all $n\in \mathbb{N}^{\ast }$, with the convention that $\underset{%
k\text{ with }a_{k}^{n}=a}{\sum }\frac{1}{2^{k}}$ (respectively $\underset{k%
\text{ with }a_{k}=a}{\sum }\frac{1}{2^{k}}$) is $0$ if there is no $k$ such
that $a_{k}^{n}=a$ (respectively $a_{k}=a$), we have

\begin{equation*}
|\alpha _{a_{1}}^{n}-\alpha _{a_{1}}|=|\sum_{k\text{ with }a_{k}^{n}=a_{1}}%
\frac{1}{2^{k}}-\sum_{k\text{ with }a_{k}=a_{1}}\frac{1}{2^{k}}|=
\end{equation*}

\begin{equation}
=|\sum_{k>1\text{ with }a_{k}^{n}=a_{1}}\frac{1}{2^{k}}-\sum_{k\text{ with }%
a_{k}=a_{1}}\frac{1}{2^{k}}|\geq \frac{1}{2}-\sum_{k>1\text{ with }%
a_{k}^{n}=a_{1}}\frac{1}{2^{k}}\text{.}  \tag{**}
\end{equation}

We claim that for every $m\in \mathbb{N}^{\ast }$, $m\geq 2$, there exists $%
n_{m}\in \mathbb{N}^{\ast }$ such that $%
a_{2}^{n}=a_{3}^{n}=...=a_{m}^{n}=a_{1}$, for every $n\in \mathbb{N}^{\ast }$%
, $n\geq n_{m}$.

Then, with the notation $\beta =za_{1}a_{1}...a_{1}...$, we have 
\begin{equation*}
d(\alpha ^{n},\beta )\leq \frac{1}{m}\text{,}
\end{equation*}
for every $n\in \mathbb{N}^{\ast }$, $n\geq n_{m}$, i.e. 
\begin{equation*}
\underset{n\rightarrow \infty }{\lim }\alpha ^{n}=\beta
\end{equation*}
(in $N(A)$).

Taking into account Theorem 2.4, we infer that $\underset{n\rightarrow
\infty }{\lim }p_{p}(\alpha ^{n})=p_{p}(\beta )$, i.e. 
\begin{equation*}
\underset{n\rightarrow \infty }{\lim }\left\Vert x_{n}-p_{p}(\beta
)\right\Vert _{p}=0\text{.}
\end{equation*}

Since, according to the hypothesis, $\underset{n\rightarrow \infty }{\lim }%
\left\Vert x_{n}-x\right\Vert _{p}=0$, we infer that $x=p_{p}(\beta )$, i.e. 
$p_{p}(\alpha )=p_{p}(\beta )$ and therefore 
\begin{equation*}
s_{p}(\widehat{\alpha })=s_{p}(\widehat{\beta })\text{.}
\end{equation*}

As $s_{p}$ is injective (see Proposition 2.1), we deduce that $\alpha $ and $%
\beta $ are equivalent, which implies the contradiction that the sequence $%
(a_{n})_{n\in \mathbb{N}^{\ast }}$\textit{\ }is constant after some rank $q$.

Consequently P($1$) is true.

Now let us prove the above claim.

Since $\underset{n\rightarrow \infty }{\lim }\left\Vert x_{n}-x\right\Vert
_{p}=0$, for every $m\in \mathbb{N}^{\ast }$, $m\geq 2$, there exists $%
n_{m}\in \mathbb{N}^{\ast }$ such that 
\begin{equation*}
\left\Vert x_{n}-x\right\Vert _{p}<\frac{1}{2^{m}}\text{,}
\end{equation*}%
for every $n\in \mathbb{N}^{\ast }$, $n\geq n_{m}$.

Then, by using $(\ast \ast )$, we get, for every $n\in \mathbb{N}^{\ast }$, $%
n\geq n_{m}$, the following inequality 
\begin{equation*}
\frac{1}{2^{m}}>\left\Vert x_{n}-x\right\Vert _{p}\geq |\alpha
_{a_{1}}^{n}-\alpha _{a_{1}}|\geq \frac{1}{2}-\sum_{k>1\text{ with }%
a_{k}^{n}=a_{1}}\frac{1}{2^{k}}\text{,}
\end{equation*}%
i.e.%
\begin{equation*}
\frac{1}{2^{m}}>\frac{1}{2}-(\sum_{k\in \{2,3,...,m\}}\frac{1}{2^{k}}%
-\sum_{k\in \{2,3,...,m\}\text{ with }a_{k}^{n}\neq a_{1}}\frac{1}{2^{k}}%
+\sum_{k\geq m+1\text{ with }a_{k}^{n}=a_{1}}\frac{1}{2^{k}})\text{,}
\end{equation*}%
which is equivalent to%
\begin{equation*}
\frac{1}{2^{m}}>\frac{1}{2}-\frac{1}{2}+\frac{1}{2^{m}}+\sum_{k\in
\{2,3,...,m\}\text{ with }a_{k}^{n}\neq a_{1}}\frac{1}{2^{k}}-\sum_{k\geq m+1%
\text{ with }a_{k}^{n}=a_{1}}\frac{1}{2^{k}}\text{,}
\end{equation*}%
from which we get%
\begin{equation*}
\sum_{k\in \{2,3,...,m\}\text{ with }a_{k}^{n}\neq a_{1}}\frac{1}{2^{k}}%
<\sum_{k\geq m+1\text{ with }a_{k}^{n}=a_{1}}\frac{1}{2^{k}}\leq \sum_{k\geq
m+1}\frac{1}{2^{k}}=\frac{1}{2^{m}}\text{.}
\end{equation*}

If there exists $s\in \{2,3,...,m\}$ such that $a_{s}^{n}\neq a_{1}$, then
by the above inequality, we obtain 
\begin{equation*}
\frac{1}{2^{s}}\leq \sum_{k\in \{2,3,...,m\}\text{ with }a_{k}^{n}\neq a_{1}}%
\frac{1}{2^{k}}<\frac{1}{2^{m}}\text{,}
\end{equation*}%
which leads us to the contradiction $m<s$.

Consequently, the claim is proved.

Now let us prove that P($m$) implies P($m+1$).

If this is not the case, then P($m$) is true and P($m+1$) is false.
Consequently, for each $n\in \mathbb{N}^{\ast }$ there exist $n^{^{\prime
}}\in \mathbb{N}^{\ast }$, $n^{^{\prime }}>n$ and $k\in \{1,2,...,m,m+1\}$
such that 
\begin{equation*}
a_{k^{^{\prime }}}^{n^{^{\prime }}}\neq a_{k^{^{\prime }}}\text{.}
\end{equation*}

Therefore, for each $n\in \mathbb{N}^{\ast }$, $n>n_{m}$, there exists $%
n^{^{\prime }}\in \mathbb{N}^{\ast }$, $n^{^{\prime }}>n$ such that $%
a_{m+1}^{n^{^{\prime }}}\neq a_{m+1}$ and hence, by passing to a
subsequence, we can suppose that 
\begin{equation*}
a_{m+1}^{n}\neq a_{m+1}\text{,}
\end{equation*}%
for every $n\in \mathbb{N}^{\ast }$, $n>n_{m}$.

If $\{n\in \mathbb{N}^{\ast }\mid a_{m+1}^{n}\neq z\}$ is infinite, by
passing to a subsequence, we can suppose that $a_{m+1}^{n}\neq z$, for every 
$n\in \mathbb{N}^{\ast }$.

We claim that the sequence $(a_{n})_{n\in \mathbb{N}^{\ast }}$ is constant
after the rank $m+1$.

This contradicts the fact that for every $q\in \mathbb{N}^{\ast }$ this
sequence\ is not constant after the rank $q$.

The claim is true since otherwise there exist $s,t,l,n\in \mathbb{N}^{\ast }$%
, such that $m+1<s<t<l$, $n_{m}<n$, $a_{s}\neq a_{t}$ and 
\begin{equation*}
\left\Vert x_{n}-x\right\Vert _{p}<\frac{1}{2^{q}}\text{,}
\end{equation*}%
for all $q\in \{m+1,m+2,...,l\}$.

If for a given $q\in \{m+1,m+2,...,l\}$ we have 
\begin{equation*}
a_{q}\neq a_{m+1}^{n}\text{,}
\end{equation*}%
then, since 
\begin{equation*}
\sum_{k>m+1\text{ with }a_{k}=a_{m+1}^{n}}\frac{1}{2^{k}}=\sum_{k>m+1\text{
with }a_{k}=a_{m+1}^{n}\text{ and }k\neq q}\frac{1}{2^{k}}<\frac{1}{2^{m+1}}-%
\frac{1}{2^{q}}\text{,}
\end{equation*}%
we get, the following contradiction%
\begin{equation*}
\frac{1}{2^{q}}<\frac{1}{2^{m+1}}-\sum_{k>m+1\text{ with }a_{k}=a_{m+1}^{n}}%
\frac{1}{2^{k}}\leq
\end{equation*}%
\begin{equation*}
\leq \left\vert \sum_{k>m\text{ with }a_{k}^{n}=a_{m+1}^{n}}\frac{1}{2^{k}}%
-\sum_{k>m+1\text{ with }a_{k}=a_{m+1}^{n}}\frac{1}{2^{k}}\right\vert =
\end{equation*}%
\begin{equation*}
=\left\vert \sum_{k\text{ with }a_{k}^{n}=a_{m+1}^{n}}\frac{1}{2^{k}}-\sum_{k%
\text{ with }a_{k}=a_{m+1}^{n}}\frac{1}{2^{k}}\right\vert =|\alpha
_{a_{m+1}^{n}}^{n}-\alpha _{a_{m+1}^{n}}|\leq \left\Vert x_{n}-x\right\Vert
_{p}<\frac{1}{2^{q}}\text{.}
\end{equation*}

Consequently we infer that 
\begin{equation*}
a_{m+1}=a_{m+2}=...=a_{s}=...=a_{t}=...=a_{l}=a_{m+1}^{n}
\end{equation*}%
which contradicts the fact that 
\begin{equation*}
a_{s}\neq a_{t}\text{.}
\end{equation*}

If $\{n\in \mathbb{N}^{\ast }\mid a_{m+1}^{n}\neq z\}$ is finite, by passing
to a subsequence, we can suppose that $a_{m+1}^{n}=z$, for every $n\in 
\mathbb{N}^{\ast }$.

As in the proof of the validity of P($1$) we obtain that $%
a_{1}a_{2}...a_{m}za_{m+1}a_{m+1}...a_{m+1}...$ and $\alpha $, where $%
a_{m+1}\neq z$, are equivalent, which implies the contradiction that the
sequence $(a_{n})_{n\in \mathbb{N}^{\ast }}$\textit{\ }is constant after
some rank $q$.

Consequently P($m+1$) is true. $\square $

\bigskip

\textbf{Proposition 3.4}. \textit{Let }$(x_{n})_{n\in \mathbb{N}^{\ast }}$ 
\textit{be a sequence of elements from }$\omega _{p}^{A}$\textit{, where} $%
x_{n}=p_{p}(\alpha ^{n})=(\alpha _{a}^{n})_{a\in A^{\prime }}$\textit{,} 
\textit{with }$\alpha ^{n}=a_{1}^{n}a_{2}^{n}...\in N(A)$,\textit{\ and }$%
x=p_{p}(\alpha )=(\alpha _{a})_{a\in A^{\prime }}\in \omega _{p}^{A}$\textit{%
, such that there exist }$n_{0}\in \mathbb{N}^{\ast }$\textit{,} $n_{0}\geq
2 $ \textit{and} $a\neq b$ \textit{such that }$\alpha
=a_{1}a_{2}...a_{n_{0}-1}abbb...\in N(A)$.

\textit{If}

\begin{equation*}
\underset{n\rightarrow \infty }{\lim }||x_{n}-x||_{p}=0\text{,}
\end{equation*}%
\textit{then for every }$m\in \mathbb{N}^{\ast }$,\textit{\ }$m>n_{0}+1$%
\textit{, there exists }$l_{m}\in \mathbb{N}^{\ast }$\textit{\ such that for
every }$l\in \mathbb{N}^{\ast }$\textit{,} $l\geq l_{m}$\textit{\ one of the
following situations is valid:}

\textit{i) }%
\begin{equation*}
a_{1}^{l}=a_{1}\text{, }a_{2}^{l}=a_{2}\text{, }...\text{, }%
a_{n_{0}-1}^{l}=a_{n_{0}-1}\text{,}
\end{equation*}%
\begin{equation*}
a_{n_{0}}^{l}=a\text{,}
\end{equation*}%
\textit{and} 
\begin{equation*}
a_{n_{0}+1}^{l}=a_{n_{0}+2}^{l}=...=a_{m}^{l}=b\text{.}
\end{equation*}

ii)\textit{\ }%
\begin{equation*}
a_{1}^{l}=a_{1}\text{, }a_{2}^{l}=a_{2}\text{, }...\text{, }%
a_{n_{0}-1}^{l}=a_{n_{0}-1}\text{,}
\end{equation*}%
\begin{equation*}
a_{n_{0}}^{l}=b\text{,}
\end{equation*}%
\textit{and} 
\begin{equation*}
a_{n_{0}+1}^{l}=a_{n_{0}+2}^{l}=...=a_{m}^{l}=a\text{.}
\end{equation*}

\textit{Proof}. As in the previous proposition one can prove that there
exists an $l_{0}\in \mathbb{N}^{\ast }$ such that 
\begin{equation*}
a_{1}^{l}=a_{1}\text{, }a_{2}^{l}=a_{2}\text{, }...\text{, }%
a_{n_{0}-1}^{l}=a_{n_{0}-1}
\end{equation*}%
for every $l\in \mathbb{N}^{\ast }$, $l\geq l_{0}$.

Let $l_{m}$ be such that $l_{m}\geq \max \{l_{0},n_{0}\}$ and 
\begin{equation*}
||x_{n}-x||_{p}<2^{-\frac{1}{q}}\frac{1}{2^{m+1}}\text{,}
\end{equation*}%
for every $n\in \mathbb{N}^{\ast }$, $n\geq l_{m}$, where $\frac{1}{p}+\frac{%
1}{q}=1$.

Let us suppose that $z\notin \{a,b\}$. The cases $z=a\neq b$ and $z=b\neq a$
could be treated in a similar way. The difference between these two cases is
the same as the difference between the two cases that occured when we
checked the validity of $P(1)$ from the proof of Proposition 3.3.

For $n$ as above, using H\"{o}lder inequality, with the convention made in
the previous proposition, we have

\begin{equation*}
||x_{n}-x||_{p}=(\sum_{A^{\prime }}|\alpha _{a}^{n}-\alpha _{a}|^{p})^{\frac{%
1}{p}}\geq (|\alpha _{a}^{n}-\alpha _{a}|^{p}+|\alpha _{b}^{n}-\alpha
_{b}|^{p})^{\frac{1}{p}}=
\end{equation*}

\begin{equation*}
=(|\sum_{k\text{ with }a_{k}^{n}=a}\frac{1}{2^{k}}-\sum_{k\text{ with }%
a_{k}=a}\frac{1}{2^{k}}|^{p}+|\sum_{k\text{ with }a_{k}^{n}=b}\frac{1}{2^{k}}%
-\sum_{k\text{ with }a_{k}=b}\frac{1}{2^{k}}|^{p})^{\frac{1}{p}}=
\end{equation*}

\begin{equation*}
=(|\sum_{k\geq n_{0}\text{ with }a_{k}^{n}=a}\frac{1}{2^{k}}-\sum_{k\geq
n_{0}\text{ with }a_{k}=a}\frac{1}{2^{k}}|^{p}+|\sum_{k\geq n_{0}\text{ with 
}a_{k}^{n}=b}\frac{1}{2^{k}}-\sum_{k\geq n_{0}\text{ with }a_{k}=b}\frac{1}{%
2^{k}}|^{p})^{\frac{1}{p}}=
\end{equation*}

\begin{equation*}
=(|\frac{1}{2^{n_{0}}}-\sum_{k\geq n_{0}\text{ with }a_{k}^{n}=a}\frac{1}{%
2^{k}}|^{p}+|\frac{1}{2^{n_{0}}}-\sum_{k\geq n_{0}\text{ with }a_{k}^{n}=b}%
\frac{1}{2^{k}}|^{p})^{\frac{1}{p}}=
\end{equation*}

\begin{equation*}
=2^{-\frac{1}{q}}\cdot 2^{\frac{1}{q}}(|\frac{1}{2^{n_{0}}}-\sum_{k\geq n_{0}%
\text{ with }a_{k}^{n}=a}\frac{1}{2^{k}}|^{p}+|\frac{1}{2^{n_{0}}}%
-\sum_{k\geq n_{0}\text{ with }a_{k}^{n}=b}\frac{1}{2^{k}}|^{p})^{\frac{1}{p}%
}\geq
\end{equation*}

\begin{equation*}
\geq 2^{-\frac{1}{q}}(|\frac{1}{2^{n_{0}}}-\sum_{k\geq n_{0}\text{ with }%
a_{k}^{n}=a}\frac{1}{2^{k}}|+|\frac{1}{2^{n_{0}}}-\sum_{k\geq n_{0}\text{
with }a_{k}^{n}=b}\frac{1}{2^{k}}|)\text{.}
\end{equation*}

We claim that 
\begin{equation*}
a_{n_{0}}^{l}\in \{a,b\}\text{,}
\end{equation*}%
for all $l\in \mathbb{N}^{\ast }$, $l\geq l_{m}$.

Indeed, if this is not the case, then, according to the above inequality, we
have%
\begin{equation*}
2^{-\frac{1}{q}}\frac{1}{2^{m+1}}>||x_{l}-x||_{p}>
\end{equation*}%
\begin{equation*}
\geq 2^{-\frac{1}{q}}(|\frac{1}{2^{n_{0}}}-\sum_{k>n_{0}\text{ with }%
a_{k}^{n}=a}\frac{1}{2^{k}}|+|\frac{1}{2^{n_{0}}}-\sum_{k>n_{0}\text{ with }%
a_{k}^{n}=b}\frac{1}{2^{k}}|)=
\end{equation*}%
\begin{equation*}
=2^{-\frac{1}{q}}(\frac{1}{2^{n_{0}-1}}-\sum_{k>n_{0}\text{ with }a_{k}^{n}=a%
\text{ or }a_{k}^{n}=b}\frac{1}{2^{k}})\geq 2^{-\frac{1}{q}}(\frac{1}{%
2^{n_{0}-1}}-\frac{1}{2^{n_{0}}})=2^{-\frac{1}{q}}\frac{1}{2^{n_{0}}}\text{,}
\end{equation*}%
so we get the contradiction%
\begin{equation*}
n_{0}\geq m+1\text{.}
\end{equation*}

If 
\begin{equation*}
a_{n_{0}}^{l}=a\text{,}
\end{equation*}%
we claim that\textit{\ }%
\begin{equation*}
a_{n_{0}+1}^{l}=a_{n_{0}+2}^{l}=...=a_{m}^{l}=b\text{.}
\end{equation*}

If this is not the case, there exists $s\in \{n_{0}+1,n_{0}+2,...,m\}$ such
that 
\begin{equation*}
a_{s}^{l}\neq b\text{.}
\end{equation*}

Then 
\begin{equation*}
\frac{1}{2^{m+1}}>\left\Vert x_{l}-x\right\Vert _{p}=(\sum_{A^{\prime
}}|\alpha _{a}^{l}-\alpha _{a}|^{p})^{\frac{1}{p}}\geq |\alpha
_{b}^{l}-\alpha _{b}|=
\end{equation*}%
\begin{equation*}
=|\sum_{k\text{ with }a_{k}^{l}=b}\frac{1}{2^{k}}-\sum_{k\text{ with }%
a_{k}=b}\frac{1}{2^{k}}|=|\sum_{k\geq n_{0}\text{ with }a_{k}^{l}=b}\frac{1}{%
2^{k}}-\sum_{k\geq n_{0}\text{ with }a_{k}=b}\frac{1}{2^{k}}|\geq
\end{equation*}%
\begin{equation*}
\geq \frac{1}{2^{n_{0}}}-\sum_{k>n_{0}\text{ with }a_{k}^{l}=b}\frac{1}{2^{k}%
}\geq \frac{1}{2^{n_{0}}}-(\frac{1}{2^{n_{0}}}-\frac{1}{2^{s}})=\frac{1}{%
2^{s}}
\end{equation*}%
and we obtain the contradiction%
\begin{equation*}
s\geq m+1\text{.}
\end{equation*}

If 
\begin{equation*}
a_{n_{0}}^{l}=b\text{,}
\end{equation*}%
in a similar manner we can prove that\textit{\ }%
\begin{equation*}
a_{n_{0}+1}^{l}=a_{n_{0}+2}^{l}=...=a_{m}^{l}=a\text{. }\square
\end{equation*}

\bigskip

\textbf{Theorem 3.1}. \textit{Let }$(x_{n})_{n\in \mathbb{N}^{\ast }}$ 
\textit{be a sequence of elements from }$\omega _{p}^{A}$\textit{, where} $%
x_{n}=p_{p}(\alpha ^{n})=(\alpha _{a}^{n})_{a\in A^{\prime }}$\textit{,} 
\textit{with }$\alpha ^{n}=a_{1}^{n}a_{2}^{n}...\in N(A)$,\textit{\ and }$%
x=p_{p}(\alpha )=(\alpha _{a})_{a\in A^{\prime }}\in \omega _{p}^{A}$\textit{%
, where }$\alpha =a_{1}a_{2}...\in N(A)$\textit{.}

\textit{i)} \textit{If for every }$q\in \mathbb{N}^{\ast }$ \textit{the
sequence }$(a_{n})_{n\in \mathbb{N}^{\ast }}$\textit{\ is not constant after
the rank }$q$\textit{, then the following assertions are equivalent:}

\textit{a)}

\begin{equation*}
\underset{n\rightarrow \infty }{\lim }||x_{n}-x||_{p}=0\text{;}
\end{equation*}

\textit{b)}%
\begin{equation*}
\underset{n\rightarrow \infty }{\lim }\alpha ^{n}=\alpha \text{,}
\end{equation*}%
\textit{in }$N(A)$\textit{;}

ii) \textit{If there exist }$n_{0}\in \mathbb{N}^{\ast }$\textit{,} $%
n_{0}\geq 2$ \textit{and} $a\neq b$ \textit{such that }$\alpha
=a_{1}a_{2}...a_{n_{0}-1}abbb....\in N(A)$, \textit{the following assertions
are equivalent:}

\textit{a)}

\begin{equation*}
\underset{n\rightarrow \infty }{\lim }||x_{n}-x||_{p}=0\text{;}
\end{equation*}

\textit{b)} \textit{for every }$m\in \mathbb{N}^{\ast }$,\textit{\ }$%
m>n_{0}+1$\textit{, there exists }$l_{m}\in \mathbb{N}^{\ast }$\textit{\
such that for every }$l\in \mathbb{N}^{\ast }$\textit{,} $l\geq l_{m}$ 
\textit{one of the following situations is valid:}

$\qquad \alpha $\textit{) }%
\begin{equation*}
a_{1}^{l}=a_{1}\text{, }a_{2}^{l}=a_{2}\text{, }...\text{, }%
a_{n_{0}-1}^{l}=a_{n_{0}-1}\text{,}
\end{equation*}%
\begin{equation*}
a_{n_{0}}^{l}=a\text{,}
\end{equation*}%
\textit{and} 
\begin{equation*}
a_{n_{0}+1}^{l}=a_{n_{0}+2}^{l}=...=a_{m}^{l}=b\text{;}
\end{equation*}

$\qquad \beta $)\textit{\ }%
\begin{equation*}
a_{1}^{l}=a_{1}\text{, }a_{2}^{l}=a_{2}\text{, }...\text{, }%
a_{n_{0}-1}^{l}=a_{n_{0}-1}\text{,}
\end{equation*}%
\begin{equation*}
a_{n_{0}}^{l}=b\text{,}
\end{equation*}%
\textit{and} 
\begin{equation*}
a_{n_{0}+1}^{l}=a_{n_{0}+2}^{l}=...=a_{m}^{l}=a\text{.}
\end{equation*}

\textit{Proof}.

i) $a)\Rightarrow b)$ results from Proposition 3.3.

$b)\Rightarrow a)$ results from Proposition 3.2.

ii) $a)\Rightarrow b)$ results from Proposition 3.4.

For $b)\Rightarrow a)$ we divide the sequence $(\alpha ^{n})_{n}$ in (at
most) two subsequences such that the first one is convergent to $\alpha
=a_{1}a_{2}...a_{n_{0}-1}abb...b...$ in $N(A)$ and the second one is
convergent to $\beta =a_{1}a_{2}...a_{n_{0}-1}baa...a...$ in $N(A)$. Then we
apply Proposition 3.2 for these subsequences of $(\alpha ^{n})_{n}$ and we
take into account the equality $p_{p}(\alpha )=s_{p}(\widehat{\alpha }%
)=s_{p}(\widehat{\beta })=p_{p}(\beta )$ which is valid since $\alpha $ and $%
\beta $ are equivalent. $\square $

\bigskip

\textbf{Remark 3.1}. \textit{Let us note that, according to the results from
this section, the topological structure of }$\omega _{p}^{A}$\textit{\ is
independent of }$p$\textit{.}

\bigskip

\begin{center}
REFERENCES

\bigskip
\end{center}

[1] S.L. Lipscomb and J.C. Perry, \textit{Lipscomb's }$L(A)$\textit{\ space
fractalized in Hilbert's }$l^{2}(A)$\textit{\ space}, Proc. Amer. Math. Soc. 
\textbf{115} (1992), 1157-1165.

[2] S.L. Lipscomb, \textit{On imbedding finite-dimensional metric spaces},
Trans. Amer. Math. Soc. \textbf{211} (1975), 143-160.

[3] S. L. Lipscomb, \textit{Fractals and universal spaces in dimension theory%
}, Springer Verlag, 2009.

[4] R. Miculescu and A. Mihail, \textit{Lipscomb space }$\omega ^{A}$\textit{%
\ is the attractor of an infinite IFS containing affine transformations of }$%
l^{2}(A)$, Proc. Amer. Math. Soc. \textbf{136} (2008), no. 2, 587--592.

[5] A. Mihail and R. Miculescu, \textit{The shift space for an infinite
iterated function system}, Math. Rep. Bucur \textbf{11 (61)} (2009), 21-32.

[6] A. Mihail and R. Miculescu, \textit{Lipscomb's }$L(A)$\textit{\ space
fractalized in Hilbert's }$l^{p}(A)$, Mediterr. J. Math. DOI
10.1007/s00009-003-0000, in press.

[7] G. N\"{o}beling, \textit{\"{U}ber eine }$n$\textit{-dimensionale
Universalmenge in} $\mathbb{R}_{2n+1}$, Math. Ann. \textbf{104} (1931),
71-80.

[8] J.C. Perry, \textit{Lipscomb's universal space is the attractor of an
infinite iterated function system}, Proc. Amer. Math. Soc. \textbf{124}
(1996), 2479-2489.

\bigskip

{\small DEPARTMENT\ OF\ MATHEMATICS, BUCHAREST\ UNIVERSITY, BUCHAREST,
ACADEMIEI\ STREET, No. 14, ROMANIA}

{\small E-mail addresses: miculesc@yahoo.com; mihail\_alex@yahoo.com.}

\end{document}